\documentclass[12pt, english]{article}
\usepackage[all]{xy}
\usepackage{varioref,babel}
\usepackage{enumerate,a4,amsfonts,theorem,amsmath,amsfonts,amssymb,amscd}
\usepackage{amsxtra,eucal}
\usepackage[latin1]{inputenc}
\usepackage{tabularx}

\begin{document}
\pagestyle{plain}
\parindent0mm

\title{Sesquilinear forms over rings with involution}
\author{Eva Bayer-Fluckiger and Daniel Arnold Moldovan\footnote{Partially supported by the
Swiss National Science Fundation, grant 200020-109174/1} \\ \'Ecole Polytechnique
F\'ed\'erale de
Lausanne, Switzerland}
\noindent
\maketitle

\textbf{Abstract }  Many classical results concerning quadratic forms have been extended to hermitian forms over algebras with involution. However, not much is known in the case
of sesquilinear forms without any symmetry property. The present paper will establish a Witt cancellation result, an analogue of Springer's theorem, as well as some local-global
and finiteness results in this context.

\bigskip

 \medskip \\


\medskip
\bigskip

{\bf \par Introduction}
\bigskip \\
Many classical results on quadratic forms
have been extended to hermitian forms over some rings with
involution  (see for instance [9], [14]). However, not much is known about sesquilinear
forms without any symmetry property. The aim of this paper is to 
show that the category of sesquilinear forms over rings with involution
is equivalent to the category of unimodular hermitian forms over a suitable additive category (cf. \S 1 --  \S 4),
and then to 
apply the theory of Quebbemann, Scharlau, Scharlau and Schulte (cf. [11], [12]) 
as well as known results concerning hermitian forms.  
\smallskip \\

This enables us to prove Witt's cancellation theorem (see theorem 6.1) and an analogue of Springer's theorem (see theorem 7.1)
for  sesquilinear forms over finite-dimensional algebras with involution over fields of characteristic different from $2$. 
We also obtain local-global results (see theorem 9.1 and its corollary) as well as finiteness theorems.

\bigskip
\bigskip

\newpage

{\bf \S 1. Sesquilinear forms over rings with involution}
\medskip \\

Let $A$ be a ring. An {\it involution} on $A$ is by definition an additive map  $\sigma : A \to A$ 
such that $\sigma(ab) = \sigma(b) \sigma(a)$ for all $a,b \in A$ and $\sigma^2$ is the identity. 
Let $V$ be a right $A$--module of finite type. 
A {\it sesquilinear form} over $(A,\sigma)$ is a biadditive map
$s : V \times V \to A$ satisfying the condition $s(xa,yb) = \sigma(a) s(x,y) b$ for all $x,y \in V$ and all $a, b
 \in A$. The {\it orthogonal sum} of two sesquilinear forms $(V,s)$ and $(V',s')$
is by definition the form $(V \oplus V', s \oplus s')$ defined by 
$$(s \oplus s')(x \oplus x', y \oplus y') =
s(x,y) + s'(x',y')$$
for all $x, y \in V$ and $x', y' \in V'$. Two sesquilinear forms $(V,s)$ and $(V', s')$ are called \textit{isometric} if
there exists an isomorphism of $A$-modules $f: V \tilde{\rightarrow} V'$ such that $s'(f(x), f(y))=s(x,y)$ for all $x, y \in V$.
\medskip \\
Let $V^*$ be the additive group ${\rm Hom}_A(V,A)$ with the  right $A$--module structure given by
$(f \cdot a)(x) = \sigma(a) f(x)$ for all $a \in A$, $f \in V^*$ and $x \in V$. We say that $V$ is
{\it reflexive} if the homomorphism of right $A$-modules $e_V : V \to V^{**}$ defined by $e_V(x)(f) = \sigma(f(x))$ for all $x \in V$ and $f \in V^*$
is bijective. 
\medskip \\
A sesquilinear form $(V,s)$ induces two homomorphisms of right $A$-modules $V \to V^*$, called its {\it left}, respectively  
{\it right adjoint}, namely $s_l : V \to V^*$ defined by $s_l(x)(y) = s(x,y)$ and
$s_r : V \to V^*$ given by $s_r(x)(y) = \sigma(s(y,x))$ for all $x, y \in V$. Note that these
are related by $s_r = s_{\ell}^* e_V$. Note also that the data $s$ and $s_r$ (as well
as $s$ and $s_{\ell}$) are equivalent. The sesquilinear form $(V,s)$ is said to be {\it unimodular} if and only if $s_{\ell}$ is
bijective (equivalently, $s_r$ is bijective). 
\medskip \\
If $(V,s)$ and $(V',s')$ are two sesquilinear forms, then an isomorphism $f : V \to V'$ defines an isometry between $(V,s)$ and $(V',s')$ if and only if we have $s_{\ell} = f^* s_{\ell}' f$.
\medskip 
Let $\mathcal{R}$ be the category of reflexive right $A$-modules.  The morphisms
of this category are homomorphisms of $A$--modules. \medskip \\
Let us denote by ${\cal S}_{\mathcal{R}}(A,\sigma)$ the category of sesquilinear forms over $(A,\sigma)$ 
defined on objects of $\mathcal{R}$. The morphisms of this category are isometries.

\bigskip
\bigskip
\bigskip

{\bf \S 2. Hermitian categories} 
\medskip \\

The aim of this section is to recall the notion of hermitian forms in additive categories
as defined in [11], [12] (see also [9], [14]). 
\smallskip \\
Let $\cal C$ be an additive category. Let
${}^* : {\cal C} \to {\cal C}$ be a duality functor, i.e. an additive contravariant functor with
a natural isomorphism $(E_C) _{C \in {\cal C}} : {\rm id} \to {}^{**}$ such that $E^*_C E_{C^*} 
= {\rm id}_{C^*}$ for all $C \in {\cal C}$. A {\it hermitian form} in the category ${\cal C}$ is a pair $(C, h)$, where $C$ is an object of ${\cal C}$
and $h: C \rightarrow C^*$ such that $h=h^* E_C$. The hermitian form is said to
be {\it unimodular} if $h$ is an isomorphism. Orthogonal
sums are defined in the obvious way. 
\medskip \\
Let $(C,h)$ and $(C',h')$ be two hermitian forms in ${\cal C}$. We say that they are {\it isometric} if there exists an isomorphism $f : C \tilde{\rightarrow} C'$ in the category $\cal C$
such that $h = f^* h' f$. 
\medskip \\
We denote by $\mathcal{H}(\cal C)$ the category of unimodular hermitian forms in the category 
$\cal C$. The morphisms are isometries.
\medskip \\
We observe that if we take $\mathcal{C}=\mathcal{R}$ then the above notion coincides with the 
notion of hermitian form over $(A, \sigma)$ defined on objects of $\mathcal{R}$. 

\bigskip
\bigskip

{\bf \S 3. An additive category}
\medskip \\
Let $(A, \sigma)$ be a ring with involution and $\mathcal{R}$ be the category of reflexive right $A$-modules. Let us consider the category $\mathcal{M}_{\mathcal{R}}$ with 
objects of the form $(V,W,f_1,f_2)$, where $V$ and $W$ are objects of $\mathcal{R}$ of finite type and 
$f_1 :  V \to W$, $f_2 : V \to W$ are homomorphisms of $A$--modules.  A {\it morphism}  $(V,W,f_1,f_2) \to (V',W',f_1',f_2')$ in  $\mathcal{M}_{\mathcal{R}}$ is a pair $(\phi,\psi)$ of  
$A$--linear homomorphisms $\phi : V \to V'$ and $\psi : W \to W'$ such that $f'_1\phi = \psi f_1$ and $f_2'\phi = \psi f_2$. The dual of $(\phi,\psi)$ is $(\psi^*,\phi^*)$. By defining direct 
sums in the obvious way we see that $\mathcal{M}_{\mathcal{R}}$ is an additive category. It is called the \textit{category of double 
arrows between objects of $\mathcal{R}$}.
\medskip \\
Let $(W^*,V^*,f_2^*,f_1^*)$ be the dual of $(V,W,f_1,f_2)$ and set $E_{(V,W,f_1,f_2)} = (e_V,e_W)$. This defines a 
duality on the category $\mathcal{M}_{\mathcal{R}}$.  We denote by $\mathcal{H}(\mathcal{M}_{\mathcal{R}})$ the category of unimodular hermitian forms in the additive category
$\mathcal{M}_{\mathcal{R}}$, as defined in the previous section. 
\bigskip
\bigskip \\
{\bf \S 4. An equivalence of categories}
\medskip \\
Let $(A, \sigma)$ be a ring with involution and $\mathcal{R}$ be the category of reflexive right $A$-modules. The aim of this section is to prove that the category of sesquilinear forms ${\cal S}_{\mathcal{R}}(A,\sigma)$
is equivalent to the hermitian category $\mathcal{H}(\mathcal{M}_{\mathcal{R}})$ defined in \S 3. We will use the methods of [3]. 
\medskip \\
We define a functor $\cal F$ from the category $\mathcal{S}_{\mathcal{R}}(A,\sigma)$ to the 
category $\mathcal{H}(\mathcal{M}_{\mathcal{R}})$ as follows. Let $(V,s)$ be an object of $\mathcal{S}_{\mathcal{R}}(A,\sigma)$,
and let $s_{\ell} : V \to V^*$ and
$s_r : V \to V^*$ be the left, respectively the right adjoint of $(V,s)$ (cf. \S 1). Then
$(V,V^*,s_{\ell},s_r)$ is an object of $\mathcal{M}_{\mathcal{R}}$ and 
the pair $(e_V,{\rm id}_{V^*})$ defines a unimodular hermitian form on $(V,V^*,s_{\ell},s_r)$;
note that the unimodularity follows from the reflexivity of $V$. 
Set ${\cal F}(V,s) = ((V,V^*,s_{\ell},s_r), (e_V, {\rm id}_{V^*}))$. For an isomorphism $\phi: (V, s) \rightarrow (W, t)$ of sesquilinear forms, where $V, W \in \mathcal{R}$, set
${\cal F}(\phi)=(\phi, \phi^{*-1})$.
\smallskip \\

\noindent
{\bf 4.1. Theorem } {\it The functor $\cal F$ is an equivalence of categories
between $\mathcal{S}_{\mathcal{R}}(A,\sigma)$ and $\mathcal{H}(\mathcal{M}_{\mathcal{R}})$. }
\smallskip \\

\noindent
{\bf Proof.} We already know that $\cal F$ sends sesquilinear forms over $(A,\sigma)$ defined on objects of $\mathcal{R}$
to objects of  $\mathcal{H}(\mathcal{M}_{\mathcal{R}})$.  Let us check that $\cal F$ sends morphisms to morphisms.
Let $\phi : (V,s) \to (W,t)$ be an isomorphism between two objects of $\mathcal{S}_{\mathcal{R}(A, \sigma)}$. Then 
$(\phi,\phi^{* -1}) : (V,V^*,s_{\ell},s_r) \to (W,W^*,t_{\ell},t_r) $ is an isomorphism in $\mathcal{M}_{\mathcal{R}}$. Moreover, as
$\phi^{**} e_V = e_W \phi$, it is also an isomorphism in 
$\mathcal{H}(\mathcal{M}_{\mathcal{R}})$.
\medskip \\
Let us define a functor $\cal G$ from $\mathcal{H}(\mathcal{M}_{\mathcal{R}})$ to $\mathcal{S}_{\mathcal{R}}(A,\sigma)$. Let
$(M,\zeta)$ be an object of $\mathcal{H}(\mathcal{M}_{\mathcal{R}})$ with $M = (V,W,f_1,f_2)$ and
$\zeta = (\phi,\psi)$. Let us define a sesquilinear form $s : V \times V \to A$
by $s(x,y) = (\psi f_2)(x)(y)$ for all $x, y \in V$ and set ${\cal G}(V,W,f_1,f_2) = (V,s)$ (in other
words, $s_{\ell} = \psi f_2$). Then
$(V,s)$ is an object of ${\cal S}_{\mathcal{R}}(A,\sigma)$  (by definition of the right $A$--module
structure of $V^*$).  For a morphism $\lambda=(\lambda_1, \lambda_2)$ 
between two objects of $\mathcal{H}(\mathcal{M}_{\mathcal{R}})$ set ${\cal G}(\lambda)=\lambda_1$.
\medskip \\
Let us check that $\cal G$ sends morphisms to morphisms. Let 
$$\lambda=(\lambda_1, \lambda_2): ((V, W, f_1, f_2), (\phi, \psi)) \rightarrow ((V', W', f'_1, f'_2), (\phi', \psi'))$$
be a morphism between two objects of $\mathcal{H}(\mathcal{M}_{\mathcal{R}})$. Then $\lambda_1: V \rightarrow V'$ and $\lambda_2: W \rightarrow W'$ are isomorphisms such that 
$\lambda_2 f_1=f'_1 \lambda_1$ and $\lambda_2 f_2=f'_2 \lambda_1$. In addition we have $(\phi, \psi)=(\lambda_1, \lambda_2)^* (\phi', \psi') (\lambda_1, \lambda_2)$, from which we
deduce that $\phi=\lambda_2^* \phi' \lambda_1$ and $\psi=\lambda_1^* \psi' \lambda_2$. It follows that 
$$\lambda_1^* \psi' f'_2 \lambda_1=\psi \lambda_2^{-1} f'_2 \lambda_1=\psi f_2,$$
so $\lambda_1$ is a morphism from $(V, \psi f_2)$ to $(V', \psi' f'_2)$. 
\medskip \\
An easy computation (using the relationship between $s_r$ and $s_{\ell}$) shows that ${\cal F}{\cal G}$ is isomorphic to the
identity in $\mathcal{H}(\mathcal{M}_{\mathcal{R}})$. Moreover, it is clear that ${\cal G}{\cal F}$ is isomorphic to the identity in 
$\mathcal{S}_{\mathcal{R}}(A,\sigma)$. Therefore ${\cal F} : \mathcal{S}_{\mathcal{R}}(A,\sigma) \to \mathcal{H}(\mathcal{M}_{\mathcal{R}})$ is
an equivalence of categories. 

\bigskip

\bigskip

{\bf \S 5. Sesquilinear forms corresponding to a given object of $\mathcal{M}_{\mathcal{R}}$}
\medskip \\
Let $A$ be a ring with an involution $\sigma$ and $\mathcal{R}$ be the category of reflexive right $A$-modules of finite type. If $(M, s)$ is a sesquilinear form over $(A, \sigma)$ defined on an object of $\mathcal{R}$, then we denote by $q(M, s)$ the corresponding object $(M, M^*, s_l, s_r)$ of the category $\mathcal{M}_{\mathcal{R}}$.
In this section we describe the set of isometry classes of sesquilinear forms over $(A, \sigma)$ corresponding by theorem 4.1 to a given object of the category $\mathcal{M}_{\mathcal{R}}$.
\medskip \\
Let us fix a sesquilinear form $(M_0, s_0)$ over $(A, \sigma)$, where $M_0 \in \mathcal{R}$, and consider the unimodular hermitian form $\eta_0 =(e_{M_0}, {\rm id}_{M_0^*})$ defined on $Q_0=q(M_0, s_0)$.
Let $E$ be the endomorphism ring of the object $Q_0$ in $\mathcal{M}_{\mathcal{R}}$, and
let us denote by $E^{\times}$ the invertible elements of $E$. The form $\eta_0$ induces an involution $\widetilde{~}$ on $E$, defined by $\widetilde{f}=
\eta_0^{-1} f^* \eta_0$ for all $f \in E$, where $f^*$ denotes the dual of $f$ in $\mathcal{M}_{\mathcal{R}}$. 
\medskip \\
On the set $E^{+}=\{f \in E^{\times} \vert~\widetilde{f}=f\}$ we define the following equivalence relation: $f \equiv f'$ if there exists a $g \in E^{\times}$ such that 
$\widetilde{g} fg=f'$. Let $H(\widetilde{~}, E^{\times})$ denote the set of equivalence classes. It is clear that $H(\widetilde{~}, E^{\times})$ is in bijection with the set of
isometry classes of unimodular hermitian forms of rank one over $(E, \widetilde{~})$. 
\medskip \\
\textbf{5.1. Theorem } \textit{The set of isometry classes of sesquilinear forms $(M, s) \in \mathcal{S}_{\mathcal{R}}(A, \sigma)$ such that $q(M, s) \simeq Q_0$ is in bijection 
with $H(\widetilde{~}, E^{\times})$.}
\medskip \\
\textbf{Proof. } For every sesquilinear form $(M, s) \in \mathcal{S}_{\mathcal{R}}(A, \sigma)$ such that $q(M, s) \simeq Q_0$ we consider a fixed isomorphism $\varphi: q(M,h) \tilde{\rightarrow} Q_0$. Through this 
isomorphism the unimodular hermitian form $(e_M, {\rm id}_{M^*})$ on $q(M, s)$ induces a hermitian form $\eta_M=\varphi^{*-1} (e_M, id_{M^*}) \varphi^{-1}$ on $Q_0$. Moreover,
$\eta_M$ is a unimodular hermitian form~: this is easy to check, noting that if $\phi : M \to M_0$
is an isomorphism, then $\phi^{**}e_M = e_{M_0} \phi$. 
We define the following map:
$$\{[(M,s)] ~\vert~(M,s) \in \mathcal{S}_{\mathcal{R}}(A, \sigma), ~q(M, h) \simeq Q_0\} \rightarrow H(\widetilde{~}, E^{\times})$$
$$[(M, s)] \mapsto\eta_0^{-1} \eta_M.$$ The hermitianity property implies
that $f = \eta_0^{-1} \eta_M$ is an element of $E^+$. It is easy to check that the class
of $f$ is independent of the choice of the isomorphism $\phi$. 
The above map is well defined since 
\begin{center}
 $(M,s) \simeq (N,t)$ if and only if $\eta_0^{-1} \eta_M \equiv \eta_0^{-1} \eta_N$.
\end{center}
This equivalence also shows its injectivity. Its surjectivity is easy to check. 
\bigskip

\bigskip

{\bf \S 6. Witt's cancellation theorem for sesquilinear forms}
\medskip \\
The aim of this section is to prove a cancellation theorem for sesquilinear forms.
Let $K$ be a field of characteristic different from $2$, let $A$ be a finite-dimensional $K$-algebra and let $\sigma$ be an involution on $A$. Let us denote by $\mathcal{R}$ the category of reflexive 
right $A$-modules which are finite dimensional $K$-vector spaces.
\medskip \\
\textbf{6.1. Theorem } \textit{Let $(M, s)$, $(M_1, s_1)$ and $(M_2, s_2)$ be sesquilinear forms over $(A, \sigma)$ defined over objects of $\mathcal{R}$ such that
$$(M_1, s_1) \oplus (M, s) \simeq (M_2, s_2) \oplus (M, s).$$ Then we have $(M_1, s_1) \simeq (M_2, s_2)$.}
\medskip \\
By the equivalence between the categories $\mathcal{S}_{\mathcal{R}}(A, \sigma)$ and $\mathcal{H}(\mathcal{M}_{\mathcal{R}})$ given by theorem 4.1, it is enough to prove that Witt's cancellation theorem holds
in the category $\mathcal{H}(\mathcal{M}_{\mathcal{R}})$. This is the purpose of proposition 6.3. For its proof we use a result of Quebbemann, Scharlau and Schulte (see [11], 3.4.1):
\medskip \\
\textbf{6.2. Theorem } 
\textit{Let $\mathcal{C}$ be a hermitian category satisfying the following conditions:
\begin{enumerate}[(i)]
 \item All idempotents of $\mathcal{C}$ split, i.e. for any object $M$ of $\mathcal{C}$ and for any idempotent $e \in {\rm End}(M)$, there exist an object $M' \in \mathcal{C}$ and morphisms
$i: M' \rightarrow M$, $j: M \rightarrow M'$ such that $ji={\rm id}_{M'}$ and $ij=e$.
\item Every object $M$ of $\mathcal{C}$ has a decomposition of the form $M \simeq N_1 \oplus ... \oplus N_r$ with $N_i \in \mathcal{C}$ indecomposable and ${\rm End}(N_i)$ a local ring in which 
$2$ is invertible.
\item For every object $M$ of $\mathcal{C}$ the ring ${\rm End}(M)$ is $J(M)$-adically complete, where $J(M)$ is the Jacobson radical of ${\rm End}(M)$.
\end{enumerate}
If $(M, h)$, $(M_1, h_1)$ and $(M_2, h_2)$ are unimodular hermitian forms in $\mathcal{C}$ such that
$$(M_1, h_1) \oplus (M, h) \simeq (M_2, h_2) \oplus (M, h),$$ 
then $(M_1, h_1) \simeq (M_2, h_2)$.
}
\medskip \\
This theorem has been originally formulated for quadratic forms in $\mathcal{C}$. But since $2$ is invertible in the endomorphism ring of every object of $\mathcal{C}$,
the categories of quadratic and unimodular hermitian forms in $\mathcal{C}$ are isomorphic. We deduce that Witt's cancellation holds for unimodular hermitian forms in $\mathcal{C}$ as well. 
\medskip \\
\textbf{6.3. Proposition } \textit{Let $(Q, \eta)$, $(Q_1, \eta_1)$ and $(Q_2, \eta_2)$ be unimodular hermitian forms in the category $\mathcal{M}_{\mathcal{R}}$ such that
$$(Q_1, \eta_1) \oplus (Q, \eta) \simeq (Q_2, \eta_2) \oplus (Q, \eta).$$
It then follows that $(Q_1, \eta_1) \simeq (Q_2, \eta_2)$. }
\medskip \\
\textbf{Proof. } As in [3], proposition 2, we check that the category $\mathcal{M}_{\mathcal{R}}$ satisfies the conditions (i), (ii) and (iii) above. It follows from theorem 6.2 
that Witt's cancellation theorem is true in the category $\mathcal{H}(\mathcal{M}_{\mathcal{R}})$.

\bigskip 
\bigskip
\bigskip

{\bf \S 7. Springer's theorem for sesquilinear forms}
\medskip \\
The classical theorem of Springer states that if two quadratic forms over a field of characteristic different from $2$ become isometric over an extension of odd degree, then
they are already isometric over the base field. In this section we prove an analogue of Springer's theorem for sesquilinear forms defined on finite-dimensional algebras with involution
over a field of characteristic $ \neq 2$. \medskip \\
Let $K$ be a field of characteristic different from $2$, $A$ be a finite-dimensional $K$-algebra and $\sigma$ be a $K$-linear involution. We also consider a finite extension $L$ of $K$,
the finite-dimensional $L$-algebra $A_L=A \otimes_K L$ and the $L$-linear involution $\sigma_L=\sigma \otimes {\rm id}_L$ on $A_L$. If $(M,s)$ is a sesquilinear form over $(A, \sigma)$, then
we denote by $(M_L, s_L)$ the sesquilinear form over $(A_L, \sigma_L)$ obtained by extension of scalars: $M_L=M \otimes_K L$ and $s_L(x \otimes a, y \otimes b)=s(x,y) \otimes ab$ for all $x, y \in M$ and 
$a,b \in L$. Let $\mathcal{R}$ be the category of reflexive right $A$-modules which are finite dimensional $K$-vector spaces.
\medskip \\
\textbf{7.1. Theorem } \textit{Suppose that $L/K$ is an extension of odd degree. Let $(M,s)$ and $(M', s')$ be two
sesquilinear forms over $(A, \sigma)$ defined on objects of $\mathcal{R}$. If $(M_L, s_L)$ and $(M'_L, s'_L)$ are isometric over $(A_L, \sigma_L)$, then $(M, s)$ and $(M', s')$ are isometric
over $(A, \sigma)$. }
\medskip \\
In order to prove this result we will use techniques of hermitian categories. 
We denote by $\mathcal{M}_{\mathcal{R}}$ (~$(\mathcal{M}_\mathcal{R})_L$~) the category of double arrows between objects of $\mathcal{R}$ (respectively objects of $\mathcal{R}$ tensorised with $L$ over $K$). There is an obvious notion
of scalar extension from the category $\mathcal{M}_{\mathcal{R}}$ to the category $(\mathcal{M}_\mathcal{R})_L$: if $Q=(V, W, f, g)$ is an object of $\mathcal{M}_{\mathcal{R}}$, set
$Q_L=(V_L, W_L, f_L, g_L)$, where $V_L=V \otimes_K L$, $W_L=W \otimes_K L$ and the morphisms $f_L$ and $g_L$
are defined by extending $f$, respectively $g$ to $L$. Note that if $Q$ and $Q'$ are two objects
of $\mathcal{M}_{\mathcal{R}}$ which become isomorphic over a field extension of $K$, then
$Q$ and $Q'$ are isomorphic (indeed, we are dealing with a linear problem). 
\medskip \\
\textbf{Proof of theorem 7.1. } The proof uses the equivalence of categories proven in \S 4. Consider the objects
$q(M, s)$ and $q(M', s')$ of the category $\mathcal{M}_{\mathcal{R}}$ and the objects $q(M_L, s_L)$ and $q(M'_L, s'_L)$ of the category $(\mathcal{M}_{\mathcal{R}})_L$. By hypothesis 
$(M_L, s_L)$ and $(M'_L, s'_L)$ are isometric, so by theorem 4.1, the objects $q(M_L, s_L) \simeq q(M,s)_L$ and $q(M'_L, s'_L) \simeq q(M',s')_L$ of $(\mathcal{M}_{\mathcal{R}})_L$ are isomorphic. 
From this it follows that the objects $q(M, s)$ and $q(M', s')$ of the category $\mathcal{M}_{\mathcal{R}}$ are isomorphic. 
\medskip \\
Denote by $Q$ the object $q(M, s)$ of $\mathcal{M}_{\mathcal{R}}$ and by $E$ the endomorphism ring of $Q$ in the category $\mathcal{M}_{\mathcal{R}}$.According to theorem 5.1, the set of isometry classes of
sesquilinear forms $(M'', s'')$ over $(A, \sigma)$ such that $q(M'', s'') \simeq Q$ is in bijection with the set $H(\widetilde{~}, E^{\times})$. This last set is in turn
in bijection with the set of isometry classes of rank one unimodular hermitian forms over $(E, \widetilde{~})$. Denote by $\mu$ and $\mu'$ the rank one unimodular hermitian forms
over $(E, \widetilde{~})$ that correspond by the above bijections to $(M, s)$, respectively $(M', s')$. Since $(M_L, s_L) \simeq (M'_L, s'_L)$, the forms $\mu$ and $\mu'$
extend to isometric hermitian forms over $(E_L, \widetilde{~}_L)$, where
$E_L=E \otimes_K L$ and $\widetilde{~}_L=\widetilde{~} \otimes {\rm id}_L$. Applying Springer's theorem for unimodular hermitian forms over $(E, \widetilde{~})$
(cf. [2], th. 4.1 or  [5], th. 2.1), we deduce that $\mu \simeq \mu'$, hence $(M, s) \simeq (M', s')$. 

\bigskip
\bigskip
\bigskip

{\bf \S 8. Sesquilinear forms defined on algebras over complete discrete valu-ation rings}

\medskip

In this section we prove results similar to those of sections 6 and 7 for sesquilinear forms defined on algebras of finite rank over complete discrete valuation rings. 

\medskip 

Let $R$ be a complete discrete valuation ring and $A$ be an $R$-algebra of finite rank with an involution denoted by $\sigma$. Suppose that there exists an element $a$ in the center of $A$ 
satisfying the property $a+\sigma(a)=1$. For example this condition is fulfilled if $2$ is invertible in $R$.

\medskip

Let $\mathcal{R}$ be the category of reflexive right $A$-modules of finite type.

\medskip

\textbf{8.1. Theorem } \textit{Let $(M,s)$, $(M_1,s_1)$ and $(M_2, s_2)$ be sesquilinear forms over $(A, \sigma)$ defined on objects of $\mathcal{R}$ such that
$$(M_1, s_1) \oplus (M, s) \simeq (M_2, s_2) \oplus (M, s).$$
It follows that $(M_1, s_1) \simeq (M_2, s_2)$. }

\medskip

\textbf{Proof } Denote by $\mathcal{M}_{\mathcal{R}}$ the category of double arrows between objects of $\mathcal{R}$, cf. \S 3. Due to the equivalence of categories between the categories $\mathcal{S}_{\mathcal{R}}(A, \sigma)$ and 
$\mathcal{H}(\mathcal{M}_{\mathcal{R}})$ given by theorem 4.1, it is enough to prove that Witt's cancellation theorem holds in the category $\mathcal{H}(\mathcal{M}_{\mathcal{R}})$. 
This follows from theorem 6.2, observing that the category $\mathcal{M}_{\mathcal{R}}$ satisfies 
conditions (i), (ii), (iii) (see e.g. [3], proposition 2).

\medskip

Next we prove a result analogous to Springer's theorem in the following setting: Let $K$ be a non-dyadic local field, $\mathcal{O}_K$ be its ring of integers and 
$A$ be an $\mathcal{O}_K$-algebra of finite rank with an $\mathcal{O}_K$-linear involution $\sigma$. Let $L$ be a finite extension of $K$, $\mathcal{O}_L$ be its ring of integers and let 
$A_L=A \otimes_{\mathcal{O}_K} \mathcal{O}_L$ with the $\mathcal{O}_L$-linear involution $\sigma_L=\sigma \otimes {\rm id}_{\mathcal{O}_L}$. If $(M,s)$ is a sesquilinear form over $(A, \sigma)$, then we
denote by $(M_L, s_L)$ the sesquilinear form over $(A_L, \sigma_L)$ obtained by extension of scalars (in particular $M_L=M \otimes_{\mathcal{O}_K} \mathcal{O}_L$).

\medskip 

Let $\mathcal{R}$ be the category of reflexive right $A$-modules which are finite dimensional $K$-vector spaces.

\medskip

\textbf{8.2. Theorem } \textit{Suppose that the extension $L/K$ is of odd degree. Let $(M,s)$ and $(M',s')$ be two sesquilinear forms over $(A, \sigma)$ with $M, M' \in \mathcal{R}$. If $(M_L, s_L)$ and $(M'_L, s'_L)$ are isometric over 
$(A_L, \sigma_L)$, then $(M,s)$ and $(M', s')$ are isometric over $(A, \sigma)$. }

\medskip

Denote by $\mathcal{M}_{\mathcal{R}}$ (~$(\mathcal{M}_{\mathcal{R}})_L$~) the category of double arrows between objects of $\mathcal{R}$ (respectively objects of $\mathcal{R}$ extended to $\mathcal{O}_L$). There is an obvious notion
of scalar extension from the category $\mathcal{M}_{\mathcal{R}}$ to the category $(\mathcal{M}_{\mathcal{R}})_L$: for an object $Q=(M, N, f, g)$ of $\mathcal{M}_{\mathcal{R}}$ set 
$Q_L=(M_L, N_L, f_L, g_L)$, where $M_L=M \otimes_{\mathcal{O}_K} \mathcal{O}_L$,
$N_L=N \otimes_{\mathcal{O}_K} \mathcal{O}_L$ and $f_L, g_L$ are defined by extending $f$, respectively $g$. The proof of theorem 8.2 is analogous to the one of theorem 7.1,
applying [7], p. 4 instead of [4], th. 2.1. 

\bigskip 
\bigskip
\bigskip

{\bf \S 9. Weak Hasse principle for sesquilinear forms}

\medskip

The classical Hasse-Minkowski theorem states that if two quadratic forms defined over a global field $k$ of characteristic different from $2$ 
become isometric over all the completions of $k$, then they are already isometric over $k$. The aim of this section is to investigate this result in the case
of sesquilinear forms defined on a finite-dimensional skew field with involution. The case of bilinear forms has been treated by Waterhouse (cf. [15]), and is based
on a classification result of Riehm (cf. [13]). 

\medskip

Let $K$ be a field of characteristic different from $2$, $D$ be a finite-dimensional skew field with center $K$ and $\sigma$ be an involution on $D$. Denote by $k$
the fixed field of $\sigma$ in $K$. Then either $k=K$ (when $\sigma$ is said to be of the \textit{first kind}) or $K$ is a quadratic extension of $k$ and the restriction of 
$\sigma$ to $K$ is the non-trivial automorphism of $K$ over $k$ (in which case $\sigma$ is said to be of the \textit{second kind} or a \textit{unitary involution}). 

\medskip

Suppose that $k$ is a global field. For every prime spot $p$ of $k$, let $k_p$ be the completion of $k$ at $p$, $K_p=K \otimes_k k_p$ and $D_p=D \otimes_k k_p$. Then $D_p$ is an algebra with center
$K_p$ and we consider on it the extended involution $\sigma_p=\sigma \otimes {\rm id}_{k_p}$. Then from any sesquilinear form $(V,s)$ over $(D, \sigma)$ we obtain by extension of scalars a sesquilinear form
$(V_p, s_p)$ over $(D_p, \sigma_p)$, where $V_p=V \otimes_k k_p$ is a free right $D_p$-module of rank ${\rm dim}_D(V)$. 

\medskip

We say that the \textit{weak Hasse principle} holds for sesquilinear
forms over $(D, \sigma)$ if any two sesquilinear forms $(V,s)$ and $(V',s')$ over $(D, \sigma)$ that become isometric over all the completions of $k$ (i.e. $(V_p, s_p) \simeq (V'_p, s'_p)$ over
$(D_p, \sigma_p)$ for every prime spot $p$ of $k$) are already isometric over $(D, \sigma)$. 

\medskip

In the sequel we will be interested in determining when the weak Hasse principle holds and we will adopt the point of view of hermitian categories. 
Let $\mathcal{R}$ be the category of finite dimensional right $D$-vector spaces. We denote by $\mathcal{M}_{\mathcal{R}}$ ( $(\mathcal{M}_{\mathcal{R}})_p$ ) the category of double arrows between objects of $\mathcal{R}$ (respectively objects of $\mathcal{R}$ extended to $k_p$).
Let us fix a sesquilinear form $s_0 \in \mathcal{S}_{\mathcal{R}}(D, \sigma)$. Consider a sesquilinear form $s \in \mathcal{S}_{\mathcal{R}}(D, \sigma)$ that becomes isometric to $s_0$ over all the completions of $k$. Then for every prime spot $p$ of $k$ we have
$(s_0)_p \simeq s_p$ and by theorem 4.1 it follows that $q((s_0)_p) \simeq q(s_p)$ for all $p$. Since $q$ commutes with base change, we obtain $q(s_0)_p \simeq q(s)_p$. 
It follows that $q(s_0) \simeq q(s)$. 

\medskip 

We denote by $E$ the endomorphism ring of $q(s_0)$ in $\mathcal{M}_{\mathcal{R}}$. The unimodular hermitian form $\eta_0 =(e_{V_0}, {\rm id}_{V_0^*})$ defined on $q(s_0)$ induces an 
involution $\widetilde{~}$ on $E$. For every $p$ let $E_p=E \otimes_k k_p$, on which we consider the involution $\widetilde{~}_p=\widetilde{~} \otimes {\rm id}_{k_p}$. By theorem 5.1 it is clear that if the localisation map
$$\Phi: H(\widetilde{~}, E^{\times}) \rightarrow \prod_p H(\widetilde{~}_p, E_p^{\times})$$
is injective, then the weak Hasse principle holds for sesquilinear forms in $\mathcal{S}_{\mathcal{R}}(D, \sigma)$. Hence in the sequel we will study the injectivity of the map $\Phi$.

\medskip 

The involution $\widetilde{~}$ induces an involution on $E/{\rm rad}(E)$, still denoted by $\widetilde{~}$. Since $H(\widetilde{~}, E^{\times}) \simeq H(\widetilde{~}, (E/{\rm rad}(E))^{\times})$ (cf. [4], lemma 2.6), 
we can suppose that ${\rm rad}(E)=0$. The ring $E$ being artinian, it is semi-simple. We can thus write $E$ as a product of simple factors:
$$E \simeq \mathcal{M}_{n_1}(D_1) \times ... \times \mathcal{M}_{n_r}(D_r),$$
where for all $1 \le i \le r$, $D_i$ is a finite-dimensional skew field and $n_i \ge 1$. Without loss of generality we can suppose that the first $m$ factors are fixed by the involution $\widetilde{~}$
and that the other factors are interchanged two by two. For every $1 \le i \le m$ we denote by $\sigma_i$ the restriction of $\widetilde{~}$ to $\mathcal{M}_{n_i}(D_i)$. 
We then obtain the isomorphism
$$H(\widetilde{~}, E^{\times}) \simeq \prod_{i=1}^m H(\sigma_i, \mathcal{M}_{n_i}(D_i)^{\times}).$$
Hence the map $\Phi$ is injective if and only if for every $1 \le i \le s$ the localisation map
$$\Phi_i: H(\sigma_i, \mathcal{M}_{n_i}(D_i)^{\times}) \rightarrow \prod_p H((\sigma_i)_p, \mathcal{M}_{n_i}((D_i)_p)^{\times}),$$
defined in the obvious way, is injective ($(D_i)_p=D_i \otimes_k k_p$ and $(\sigma_i)_p=\sigma \otimes {\rm id}_{k_p}$). We thus obtain:

\medskip

\textbf{9.1. Theorem } \textit{If for all $1 \le i \le m$, either the involution $\sigma_i$ is unitary, or $D_i$ is a quaternion algebra and $\sigma_i$ is a symplectic involution, then
the map $\Phi$ is injective and hence the weak Hasse principle holds for sesquilinear forms over $(D, \sigma)$.}
\medskip \\
\textbf{Proof. } It suffices to prove that for all $1 \le i \le m$, the map $\Phi_i$ is injective. Let $1 \le i \le m$ and suppose that $\sigma_i$ is unitary. According to the Albert-Riehm-Scharlau theorem
([10], theorem 3.1, p. 31), there exists a unitary involution $\tau_i$ on $D_i$ with the same restriction to the center of $D_i$ as $\sigma_i$. Then by Morita theory 
([9], Chap. 1, theorem 9.3.5), the unimodular hermitian forms of rank one over $(\mathcal{M}_{n_i} (D_i), \sigma_i)$ (which are parametrised by the set
$H(\sigma_i, \mathcal{M}_{n_i}(D_i))$~) are in bijection with the unimodular hermitian forms of rank $n_i$ over $(D_i, \tau_i)$. An analogous statement holds for 
the unimodular hermitian forms of rank one over $(\mathcal{M}_{n_i} ((D_i)_p), (\sigma_i)_p)$. As the weak Hasse principle holds for unimodular hermitian forms over
a skew field with a unitary involution (cf. [14], theorem 10.6.1), the map $\Phi_i$ is injective.

\medskip

Consider now the case when $D_i$ is a quaternion algebra and $\sigma_i$ is a symplectic involution. Then we consider the canonical involution $\tau_i$ on $D_i$ and 
we proceed as above. The weak Hasse principle holds for unimodular hermitian forms over quaternion algebras with the canonical involution. Indeed, let $h$ and $h'$ be two unimodular hermitian 
forms over $(D_i, \tau_i)$. Then their trace forms $q_{h}$ and $q_{h'}$ are two quadratic forms defined on the center of $D_i$. Since they become isometric everywhere locally, 
they are isometric by the usual weak Hasse principle. But two unimodular hermitian forms over $(D_i, \tau_i)$ are isometric if and only if their trace forms are isometric. Hence it follows that
$h \simeq h'$. We conclude that the weak Hasse principle holds for unimodular hermitian forms over quaternion algebras with the canonical involution and thus the map $\Phi_i$ is injective.

\medskip 

\textbf{9.2. Corollary } \textit{If $\sigma$ is a unitary involution, then the weak Hasse principle holds for sesquilinear forms in $\mathcal{S}_{\mathcal{R}}(D, \sigma)$.}

\medskip

\textbf{Proof. } The statement immediately follows from theorem 9.1, since if $\sigma$ is unitary, then all the involutions $\sigma_i$ are unitary as well. 

\bigskip

\bigskip

{\bf \S 10. Finiteness results}

\bigskip
In this section we generalize the finiteness results of [4] to the setting of sesquilinear forms. For a ring $A$ we denote by $T(A)$ the $\mathbb{Z}$-torsion subgroup of $A$. If $R$ is a ring then we say that $A$ is \textit{$R$-finite}
if $A_R=A \otimes_{\mathbb{Z}} R$ is a finitely generated $R$-module and $T(A)$ is finite. \medskip \\
Let $(A, \sigma)$ be a ring with involution, $\mathcal{R}$ be the category of reflexive right $A$-modules of finite type
and $\mathcal{M}_{\mathcal{R}}$ be the category defined in \S 3. 
We recall from \S 4 that the functor
$$\mathcal{F}: \mathcal{S}_{\mathcal{R}}(A, \sigma) \rightarrow \mathcal{H}(\mathcal{M}_{\mathcal{R}})$$
$$(V, s) \mapsto ((V, V^*, s_l, s_r), (e_V, {\rm id}_{V^*}))$$
is an equivalence of categories. Fix a sesquilinear form $(V, s) \in \mathcal{S}_{\mathcal{R}}(A, \sigma)$ and denote by $q(V,s)$ the corresponding object 
$(V, V^*, s_l, s_r)$ of $\mathcal{M}_{\mathcal{R}}$ and by $E$ the endomorphism ring of $q(V, s)$ in $\mathcal{M}_{\mathcal{R}}$. 
\medskip \\
\textbf{10.1. Theorem } \textit{If there exists a non-zero integer $m$ such that ${\rm End}_A(V)$ is $\mathbb{Z}[1/m]$-finite, then $(V, s)$ contains only finitely many
isometry classes of orthogonal summands. }
\medskip \\
\textbf{Proof. } From the fact that ${\rm End}_A(V)$ is $\mathbb{Z}[1/m]$-finite it follows that $E$ is $\mathbb{Z}[1/m]$-finite and hence
$\mathbb{Q}$-finite. Hence by [4], theorem 1.1, $q(V,s)$
has only finitely many isometry classes of direct summands. \medskip \\
Let $(V',s')$ be an orthogonal summand of $(V,s)$. Then $q(V',s')$ is a direct summand of $q(V,s)$. Since there is an obvious surjection ${\rm End}(q(V,s)) \rightarrow {\rm End}(q(V',s'))$, we obtain
that ${\rm End}(q(V',s'))$ is $\mathbb{Z}[1/m]$-finite too. By the next proposition there exist only finitely many isometry classes of sesquilinear forms $(V'',s'') \in \mathcal{S}_{\mathcal{R}}(A,\sigma)$
such that $q(V'',s'') \simeq q(V',s')$. Hence the assertion of the theorem follows from the next result:
\medskip \\
\textbf{10.2. Proposition } \textit{Let $N$ be an object of $\mathcal{M}_{\mathcal{R}}$ and assume that there exists a non-zero integer $m$ such that ${\rm End}(N)$
is $\mathbb{Z}[1/m]$-finite. Then there exist only finitely many isometry classes of sesquilinear forms $(V,s) \in \mathcal{S}_{\mathcal{R}}(A, \sigma)$ 
such that $q(V,s) \simeq N$.}
\medskip \\
\textbf{Proof. } By theorem 4.1, every sesquilinear form $(V,s) \in \mathcal{S}_{\mathcal{R}}(A, \sigma)$ defines a unimodular hermitian form 
on $q(V,s)$. But by [4], theorem 1.2 the number of isometry classes of
unimodular hermitian forms on $N$ is finite. Hence the assertion follows. 
\medskip \\
We say that two sesquilinear forms over $(A, \sigma)$ are in the same \textit{genus} if they become isometric over $A \otimes_{\mathbb{Z}} \mathbb{Z}_p$ for all primes of $\mathbb{Q}$.
\medskip \\
\textbf{10.3. Theorem } \textit{If ${\rm End}_A(V)$ is $\mathbb{Q}$-finite then the genus of $(V, s)$ contains only a finite number of isometry classes of sesquilinear forms. }
\medskip \\
\textbf{Proof. } Since ${\rm End}_A(V)$ is $\mathbb{Q}$-finite, $E$ is $\mathbb{Q}$-finite too. Then by [4], theorem 3.4 the genus of $\mathcal{F}(V, s)$ contains only finitely many isometry classes of 
unimodular hermitian forms. The functor $\mathcal{F}$ induces a bijection between the genus of $(V,s)$ and the genus of $\mathcal{F}(V, s)$, hence the assertion follows.

\bigskip

\bigskip

{\bf \S 11. Bilinear forms invariant by a group action}

\medskip
Let $R$ be a commutative ring and $G$ be a finite group. Denote by $R[G]$ the group ring of $G$ over $R$. A {\it $G$--bilinear form} over $R$
is by definition a pair $(M,b)$, where $M$ is a right $R[G]$--module which is an $R$-module of finite type and 
$b : M \times M \to R$ is an $R$--bilinear form such that 
$b(xg, yg) = b(x, y)$ for all $g \in G$ and $x,y \in M$. 
We say that two $G$--bilinear forms $(M,b)$ and 
$(M',b')$ are \textit{isometric} if there exists an isomorphism of $R[G]$--modules $f : M \tilde{\rightarrow} M'$
such that $b'(f(x), f(y)) = b(x, y)$ for all $x,y \in M$. Two $G$-bilinear forms over $R$ are said to be in the same \textit{genus} if they become isometric over $R \otimes_{\mathbb{Z}} \mathbb{Z}_p$
for all primes of $\mathbb{Q}$.

\medskip

The aim of this section is to apply the results of sections 6, 7 and 10 to $G$--bilinear forms. 

\medskip

We endow the group ring $R[G]$ with the canonical involution $\sigma$, which is the $R$-linear involution such that $\sigma(g)=g^{-1}$ for every $g \in G$. 
Analogously to [8], theorem 7.1 it is straightforward to prove that there 
is a correspondence between $G$--bilinear forms over $R$ and
sesquilinear forms over $(R[G], \sigma)$. Indeed, if $(M, b)$ is a $G$--bilinear form over $R$ then the form
$$S: M \times M \rightarrow R[G], ~(x, y) \mapsto \sum_{g \in G} b(xg, y) g$$
for all $x, y \in M$ is clearly sesquilinear with respect to $\sigma$. Conversely, if $(M, S)$ is a sesquilinear form over $(R[G], \sigma)$ then $M$ has an obvious structure of $R$-module
and the form $ P \circ S: M \times M \rightarrow R$ is bilinear, where $P$ is the projection of $R[G]$ onto the $e_G$ component. 

\medskip

\textbf{11.1. Theorem } \textit{Suppose that $R$ is a field of characteristic $\not = 2$. Let
$(M_1,b_1)$, $(M_2,b_2)$ and $(M,b)$ be $G$--bilinear forms over $R$ such that
$(M_1,b_1) \oplus (M,b) \simeq (M_2,b_2) \oplus (M,b)$. Then
$(M_1,b_1) \simeq (M_2,b_2)$.}

\medskip

\textbf{Proof. } This follows from theorem 6.1 and the above correspondence.

\medskip

\textbf{11.2. Theorem } \textit{Suppose that $R$ is a field of characteristic $\not = 2$. If
two $G$--bilinear forms become isometric over an extension
of $R$ of odd degree then they are already isometric over $R$.}

\medskip

\textbf{Proof. } The result follows  from theorem 7.1 and the above correspondence.

\medskip

\textbf{11.3. Theorem } \textit{Let $(M, b)$ be a $G$--bilinear form over the commutative
ring $R$. If there exists a non-zero integer $m$ such that ${\rm End}_{R[G]}(M)$ is $\mathbb{Z}[1/m]$-finite, then
$(M, b)$ contains only finitely many isometry classes of orthogonal summands.}

\medskip

\textbf{Proof. } This follows from theorem 10.1 and the above correspondence.

\medskip

\textbf{11.4. Theorem } \textit{Let $(M, b)$ be a $G$--bilinear form over $R$. If ${\rm End}_{R[G]}(M)$ is $\mathbb{Q}$-finite, then the genus of $(M, b)$ contains 
only finitely many isometry classes of $G$--bilinear forms.}

\medskip

\textbf{Proof. } This follows from theorem 10.3 and the above correspondence.

\bigskip

\bigskip

{\bf \S 12. Generalization to systems of sesquilinear forms}

\medskip
Let $A$ be a ring with an involution $\sigma$ and $I$ be a set. 
A \textit{system of sesquilinear forms} over $(A, \sigma)$ is $(V, (s_i)_{i \in I})$, where $V$ is a reflexive right $A$-module of finite type and for all $i \in I$, 
$(V, s_i)$ is a sesquilinear form over $(A, \sigma)$. 
A morphism between two systems of sesquilinear forms $(V, (s_i)_{i \in I})$ and $(V', (s'_i)_{i \in I})$ consists of an isomorphism of $A$-modules
$f:V \tilde{\rightarrow} V'$ such that for every $i \in I$ and every $x, y \in V$ we have $s'_i(f(x), f(y))=s_i(x,y)$. 

\medskip

The results of the preceding sections can be generalized to systems of sesquilinear forms (for related results on systems of quadratic forms, see [1]). It is necessary to modify the definition of 
the category $\mathcal{M}_{\mathcal{R}}$:
its objects will be of the form $(V, W, (f_i, g_i)_{i \in I})$ and the morphisms between two such objects will be the obvious ones. To a 
system $(V, (s_i)_{i \in I})$ of sesquilinear forms over $(A, \sigma)$, where $V \in \mathcal{R}$, it will correspond the object $(V, V^*, (((s_i)_l, (s_i)_r))_{i \in I})$ of the 
category $\mathcal{M}_{\mathcal{R}}$. 

\bigskip

\bigskip

\textbf{Acknowledgements } The second named author would like to thank Emmanuel Lequeu for many interesting and useful discussions.

\newpage

{\bf Bibliography}
\bigskip \\

\noindent
[1] E. Bayer--Fluckiger, Principe de Hasse faible pour les syst\`emes de formes quadratiques, 
{\it J. reine angew. Math.} {\bf 378} (1987), 53-59. \\

\noindent
[2] E. Bayer--Fluckiger, Self--dual normal bases, {\it Indag. Math.} {\bf 92} (1989), 379--383.\\

\noindent
[3] E. Bayer--Fluckiger, L. Fainsilber, Non unimodular hermitian forms,
{\it Invent. Math.} {\bf 123} (1996), 233--240. \\

\noindent
[4] E. Bayer--Fluckiger, C. Kearton, S.M. Wilson, Hermitian forms in additive categories: finiteness results, 
{\it J. Algebra} {\bf 123} (1989), no. 2, 336--350. \\

\noindent
[5] E. Bayer--Fluckiger, H.W. Lenstra, Jr,  Forms in odd degree
extensions and self-dual normal bases,
{\it Amer. J. Math.} {\bf 112} (1990), 359-373. \\

\noindent
[6] C.W. Curtis, I. Reiner, \textit{Methods of representation theory, with applications to finite groups and orders, I}, Wiley, New York (1981). \\

\noindent
[7] L. Fainsilber, Formes hermitiennes sur des alg\`ebres sur des anneaux locaux, Publications Math\'ematiques de Besan\c con, 1994. \\

\noindent
[8] A. Fr\"{o}hlich, A.M. McEvett, Forms over rings with involution, {\it J. Algebra} {\bf 12} (1969), 79--104. \\

\noindent
[9] M. Knus, {\it Quadratic and hermitian forms over rings}, Grundlehren der math. Wiss.
{\bf 294}, Springer--Verlag (1991). \\

\noindent
[10] M. Knus, A. Merkurjev, M. Rost, J.-P. Tignol, {\it The book of involutions},
Coll. Pub. {\bf 44}, Amer. Math. Soc. (1998).  \\

\noindent
[11] H.--G. Quebbemann, W. Scharlau, M. Schulte, Quadratic and hermitian
forms in additive and abelian categories, {\it J. Algebra} {\bf 59} (1979), 264--289. \\

\noindent
[12] H.--G. Quebbemann, R. Scharlau, W. Scharlau, M. Schulte, Quadratische
Formen in additiven Kategorien, {\it Bull. Soc. Math. France} M\'emoire {\bf 48} (1976), 93--101. \\

\noindent

[13] C. Riehm, The equivalence of bilinear forms, {\it J. Algebra} {\bf 31} (1974), 45--66. \\

\noindent
[14] W. Scharlau, {\it Quadratic and Hermitian Forms}, Grundlehren der math. Wiss. {\bf 270},
Springer--Verlag (1985). \\

\noindent
[15] W.C. Waterhouse, A nonsymmetric Hasse-Minkowski theorem, {\it Amer. J. Math.}  {\bf 99} (1977), no. 4, 755--759.
\bigskip

\bigskip 
Eva Bayer-Fluckiger\\
\'Ecole Polytechnique F\'ed\'erale de Lausanne\\
EPFL-FSB-MATHGEOM-CSAG\\
Station 8, 1015 Lausanne\\
Switzerland \\ 
eva.bayer@epfl.ch
\bigskip \\
Daniel Arnold Moldovan \\
\'Ecole Polytechnique F\'ed\'erale de Lausanne\\
EPFL-FSB-MATHGEOM-CSAG\\
Station 8, 1015 Lausanne\\
Switzerland \\
danielarnold.moldovan@epfl.ch

\vfill
\eject
\end{document}